\documentclass{ifacconf}

\usepackage{xcolor}
\usepackage{enumerate}
\usepackage{amsmath, amssymb,mathrsfs}
\usepackage{natbib}            
\usepackage{graphicx}          

\usepackage{tikz,environ}
\usetikzlibrary{arrows,positioning,patterns,decorations.pathreplacing}

\tikzset{
    at xy split/.style 2 args={
        at={(#1,#2)}
    },
    a/.style={circle, draw=red},`'
    b/.style={rectangle, draw=blue}
}

\makeatletter
\newsavebox{\measure@tikzpicture}
\NewEnviron{scaletikzpicturetowidth}[1]{%
  \def\tikz@width{#1}%
  \begin{lrbox}{\measure@tikzpicture}%
  \BODY
  \end{lrbox}%
  \pgfmathparse{#1/\wd\measure@tikzpicture}%
  \BODY
}
\makeatother

\allowdisplaybreaks[4]

\begin{document}

\begin{frontmatter}
\title{Data-driven approximation for feasible regions in nonlinear model predictive control}
%
\thanks[footnoteinfo]{This work was supported by the National Science Foundation of China
	(Grant No. 61973214, 61433002, 61590924, 61573239).}
\author[First]{Yuanqiang Zhou}
\author[First]{Dewei Li}
\author[First]{Yugeng Xi}
\author[First]{Yunwen Xu}


\address[First]{Department of Automation, Shanghai Jiao Tong University, and Key Laboratory of System Control and Information Processing, Ministry of Education of China, Shanghai, 200240, China. \\e-mail: \{zhouyuanqiang, dwli, ygxi, willing419\}@sjtu.edu.cn}



\begin{abstract}                          
This paper develops a data-driven learning framework for approximating the feasible region and invariant set of a nonlinear system under the nonlinear Model Predictive Control (MPC) scheme. 
The developed approach is based on the feasibility information of a point-wise data set using low-discrepancy sequence. 
Using kernel-based Support Vector Machine (SVM) learning, we construct outer and inner approximations of the boundary of the feasible region and then, obtain the feasible region of MPC for the system. 
Furthermore, we extend our approach to the perturbed nonlinear systems using set-theoretic method. 
Finally, an illustrative numerical example is provided to show the effectiveness of the proposed approach.
\end{abstract}

\end{frontmatter}

\section{Introduction}


When it comes to MPC strategy for continuous-time systems, we periodically compute an input signal from solving an optimization problem and apply this signal to the system over a bounded horizon, see e.g. \cite{Rawlings2009model,xi2019predictive}. 
With the improving of control theory and technology, the new problems which involve large transients often require Nonlinear MPC (NMPC), see e.g. \cite{allgower2012nonlinear,grune2017nonlinear}. 
NMPC for continuous-time systems with input constraints is proposed, see e.g. \cite{chen1998quasi}, which considers a nonlinear system given by
\begin{align}\label{eq:system:chen1997}
\dot{x}(t) & = f(x(t), u(t)), \quad x(0) = \bar{x}, \quad t \in \mathbb{R}_{\ge 0}, \\
& \quad  x(t) \in \mathbb{R}^n, \quad u (t) \in \mathbb{U} \subset \mathbb{R}^m, 
\end{align}
where ${x}$ and $u$ are the state vector and the control input, respectively; $\mathbb{U}$ is the input constraint set and the mapping $f(\cdot, \cdot )$ satisfy $f(0, 0) = 0$. Besides, $f: \mathbb{R}^n \times \mathbb{U} \to \mathbb{R}^n$ is assumed to be twice continuously differentiable, which leads to unique solutions under the piece-wise right-continuous control inputs, and $\mathbb{U}$ is assumed to be compact, convex, and contains $0$ in its interior.

The quasi-infinite horizon NMPC developed in \cite{chen1998quasi} consists of designing a control invariant region $\Omega$, which is called as \emph{terminal region}, and a control Lyapunov function. 
The set $\Omega$ must be chosen such that it is invariant for \eqref{eq:system:chen1997} controlled by a local auxiliary controller. It can be determined off-line such that invariance property of $\Omega$ holds and all the constraints are satisfied. The actual implementation of the NMPC developed in \cite{chen1998quasi} would not use the local auxiliary controller, but rather guarantee recursive feasibility of the closed-loop system. As a contrast, the dual-mode NMPC developed in \cite{Michalska:1993} is actually a receding horizon implementation of a switching control mechanism; that is, it applies the optimal control input from solving a optimization problem over a finite horizon so as to drive the state to the terminal region $\Omega$; once the state enters $\Omega$, it switches the control input to the local feedback law so as to steer the state to the origin. 
We highlight that the terminal region $\Omega$ is a positively invariant set for \eqref{eq:system:chen1997} (see e.g.~\cite{bertsekas71}) and the analytic properties of positively invariant sets have been widely studied, see e.g.~\cite{blanchini:2007}. The design method of the set $\Omega$ developed in \cite{chen1998quasi} is based on the results presented in \cite{Michalska:1993}, which is calculated off-line by linearizing the system \eqref{eq:system:chen1997}. This method is quite conservative, especially for systems with high nonlinearities. The terminal region $\Omega$ is very important for stability analysis of the closed-loop NMPC and affects the feasible region of NMPC. 
Therefore, this paper considers the data-driven approximations for feasible regions and the terminal regions of the nonlinear systems. 

In this paper, we develop a data-driven learning framework for approximating feasible region and invariant set of a perturbed nonlinear systems using feasibility information of data samples. 
This work can be seen as an extension of the results developed in \cite{ong2006enlarging} and \cite{chakrabarty2016support,chakrabarty2020active}. 
Here, we consider the continuous-time nonlinear system with bounded disturbance and polyhedral state and input constraints. 
The main contributions of this paper are as follows: (i) we develop a computation method to obtain the invariant set for the nonlinear system given by \eqref{eq:system:chen1997}; (ii) we develop a data-driven learning method to obtain an inner and outer approximation of the feasible region of the NMPC using support vector machine; (iii) we extend the above two results to perturbed nonlinear system, where disturbance is bounded and the state and input space are polyhedral. 


\emph{Notation}:  
We denote $\mathbb{R}^n$ for the $n$-dimensional Euclidean space and $\mathcal{C}^1(\mathbb{R}^n)$ for the space of continuously differentiable functions on $\mathbb{R}^n$. Besides, we let $\mathbb{ R}$ and $\mathbb{ R}_{\ge 0}$ denote the real and the non-negative real numbers. For a set $\mathbb{S}$, we let $\mathrm{Card}(\mathbb{S})$ denote the cardinality of $\mathbb{S}$, $\mathrm{Vol}(\mathbb{S})$ denote the Lebesgue measure of $\mathbb{S}$, and $\mathrm{Int}(\mathbb{S})$ denote the set of interior points of $\mathbb{S}$. Furthermore, a continuous function $\alpha: \mathbb{R}_+ \rightarrow \mathbb{R}_+$ belongs to class $\mathcal{K}$ if it is continuous, $\alpha(0) = 0$, and it is strictly increasing. It belongs to class $\mathcal{K}_{\infty}$ if it belongs to class $\mathcal{K}$ and it is unbounded. Finally, given a point $x \in\mathbb{R}^{n}$ and a bounded and non-empty set $\mathcal{A}\subseteq\mathbb{R}^{n}$, the distance of  $x$ to the set $\mathcal{A}$ is defined as $|x|_{\mathcal{A}} = \mathrm{dist}(x, \mathcal{A})   \triangleq    \mathrm{inf}  \{ \|x-\xi\|_{2}: \xi \in\mathcal{A}\}$. A closed ball centralized at a given point $\xi_{0}\in\mathbb{R}^{n}$ with a radius of $\rho \ge 0$ is denoted as $ \mathbb{B}^n(\xi_{0}, \rho) \triangleq \{ \xi \in \mathbb{R}^n : \| \xi - \xi_0\| \le \rho \}$. If $\xi_0 = 0$, we let $\mathbb{B}^n(\rho)$ denote $ \mathbb{B}^n(0, \rho)$ as a shorthand.


\section{Problem formulation}
We extend \eqref{eq:system:chen1997} to a nonlinear perturbed system described as
\begin{align}
\dot{x} {(t)} & = {F} \left(x (t), u (t), w (t) \right), \quad x(0) = \bar{x},  \quad t \in \mathbb{R}_{\ge 0}, \label{eq:system:description} \\
& \quad  x(t) \in \mathbb{X} \subset \mathbb{R}^n, \quad u (t) \in \mathbb{U} \subset \mathbb{R}^m, \label{eq:system:description:2}
\end{align}
where ${x} \in \mathbb{R}^{n}$, $u\in \mathbb{R}^{m}$, and $w \in \mathbb{R}^{n}$ are the state vector, the control input, and the additive disturbance. Besides, we confine the state space $\mathbb{X}$ and the input space $\mathbb{U}$ to be polytopes. We give the following assumptions for \eqref{eq:system:description}: 

\begin{assum} \label{eq:J:l} 
\begin{description}
	\item[A1] The set $\mathbb{U} \in \mathbb{R}^{m}$ is compact, convex and contains $0 \in \mathbb{R}^{m} $ in its interior. The set $\mathbb{X}$ is a $n$-dimensional hypercube and contains $0 \in \mathbb{R}^{n}$ in its interior.
	\item[A2] The disturbance $w(\cdot) \in \mathbb{W}$, where $\mathbb{W} = \left\{w: |w| \le \bar w \right\}$. 
	\item[A3] $F: \mathbb{X} \times \mathbb{U} \times \mathbb{W} \to \mathbb{X}$ is locally Lipschitz in $x$ such that $\|F(x_1, u, w)-F(x_2, u, 0) \|  \le L_F  \|x_1 - x_2 \| + \beta(|w|_{\infty})$	for all ${x_1, x_2} \in \mathbb{X}$, ${u} \in \mathbb{U}$, and $w \in \mathbb{W}$, where $L_F$ is local Lipschitz constant and $\beta(\cdot)$ is a $\mathcal{K}$ function. 
\end{description}
\end{assum} 

The nominal model of the system \eqref{eq:system:description} is given by
\begin{align}\label{eq:system:nominal}
\dot{\hat{x}} (t) & = {F} \left(\hat x (t), u(t), 0 \right) \triangleq f(\hat x (t), u (t)), \notag \\
& \quad \qquad \qquad \hat x (0) = \bar{x},  \quad t \in \mathbb{R}_{\ge 0},
\end{align}
and the additive disturbance $w(t), t\ge 0$, in \eqref{eq:system:description} satisfies 
\begin{align}\label{eq:5}
w(t) & = \dot{x} (t) - f(x (t), u (t)). 
\end{align}
Assumption A2 indicates that $w(t)$ is bounded on $\mathbb{W}$. This kind of model uncertainty has been studied in previous papers about robustness in MPC as in \cite{zou2019event} and \cite{lu2019stochastic}. 

Given a state measurement $x_k$, a prediction horizon $T$, a state cost $\Psi(\cdot, \cdot)$, a terminal penalty $\Phi(\cdot)$, and a terminal set $\Omega$, an optimization problem of NMPC for \eqref{eq:system:description} is defined as follows: 
	\begin{subequations}\label{eq:7}
	\begin{align}\label{FHOCP}
	\hat{u}^\star = \underset{ u } { \arg\min }  ~ \mathcal{J}(x_k,  \hat{u}(\cdot)), 
	\end{align}
	with
	\begin{align}\label{Perf}
	\mathcal{J}(x_k,  \hat{u}(\cdot))  \buildrel \vartriangle \over  = \int_{0}^{T} \Psi(\hat x(\tau), \hat u(\tau) ) \mathrm{d} \tau + \Phi(\hat x(T)), 
	\end{align}
	subject to
		\begin{align}
		&  \dot{\hat x} (\tau) =   f(\hat{x} (\tau), \hat{u} (\tau)),  \quad 0 \le \tau \le T,   \label{eq:7A}  \\
		&   \hat{x} (\tau) \in \mathbb{X}, \quad \hat{u} (\tau) \in \mathbb{U}, \quad 0 \le \tau \le T,   \label{eq:7B} \\
		&  \hat{x} (0) = x_k,  \quad  \hat x (T) \in \Omega. \label{eq:7C}
		\end{align}
	\end{subequations}
For \eqref{eq:7}, we give the following two assumptions. 
\begin{assum} \label{eq:J:II}
	The stage cost $\Psi(x, u)$ in \eqref{Perf} satisfies $\Psi(0, 0) = 0$ and $\Psi(x, u) \ge \alpha_\Psi (\|x\|) + \rho_\Psi$, where $\rho_\Psi\in \mathbb{R}$ and $\alpha_\Psi(\cdot)$ is a $\mathcal{K}_{\infty}$ function. Moreover, $\Psi(x, u)$ is locally Lipschitz in $x$ in $\mathbb{X} \times \mathbb{U}$ with constant $L_\Psi$; that is, $\| \Psi(x_1, u) - \Psi(x_2, u) \| \le L_\Psi \|x_1 - x_2\|$ for all $x_1, x_2 \in \mathbb{X}$ and all $u \in \mathbb{U}$. 
\end{assum} 
\begin{assum} \label{eq:J:F}
	The terminal cost $\Phi(x)$ in \eqref{Perf} satisfies
	\begin{enumerate}
		\item $\underline{\alpha}_{\Phi} (\|x\|) \le \Phi(x) \le \overline{\alpha}_{\Phi} (\|x\|)$, where $\underline{\alpha}_{\Phi} (\cdot)$ and $\overline{\alpha}_{\Phi} (\cdot)$ are $\mathcal{K}_{\infty}$ functions;
		\item $\Phi(f(x,\kappa(x))) - \Phi(x) \le - \Psi(x, \kappa(x))$ for all $x \in \Omega$, where $\kappa(\cdot): \mathbb{ R}^n \to \mathbb{ R}^m$ is continuous and $\kappa(0) = 0$; 
		\item $\Phi(\cdot)$ is locally Lipschitz in $\Omega$ with a Lipschitz constant $L_\Phi$. 
	\end{enumerate}
\end{assum} 
Assumptions~\ref{eq:J:II} and \ref{eq:J:F} are fairly standard assumptions to guarantee feasibility and stability for nonlinear systems under MPC framework, see e.g. \cite{allgower2012nonlinear,grune2017nonlinear}. Specifically, several numerical methods to obtain $\Omega$ and the local auxiliary controller $\kappa(x)$ satisfying Assumption \ref{eq:J:F} were proposed in \cite{Michalska:1993,chen1998quasi}. The feature of $\kappa(x)$ makes it desirable to construct the feasible control input and further to prove the closed-loop stability. Furthermore, several ways to design $\Psi(\cdot, \cdot)$ and $\Phi(\cdot)$ satisfying Assumptions~\ref{eq:J:l} and \ref{eq:J:II} were proposed in \cite{Hashimoto2016self}.  In practice, $\Psi(x, u)$ and $\Phi(x)$ can be taken as quadratic functions.

Now, we assume that for \eqref{eq:7}, Assumptions~\ref{eq:J:l} and \ref{eq:J:II} hold. Then we introduce the following definitions. 
\begin{defn}[Feasible State]  \label{Feasible:state}
	A sampling state $x_k \in \mathbb{X}$ is feasible if there exists a admissible control function $\kappa(\cdot) \in \mathbb{U}$, $t \ge 0$, satisfying all the conditions  \eqref{eq:7A}--\eqref{eq:7C} except the minimization of \eqref{Perf}. 
\end{defn}
Using Definition~\eqref{Feasible:state}, we define the feasible region $\mathbb{F}$: 
\begin{equation} \label{definition:F}
\mathbb{F} \buildrel \vartriangle \over  = \big\{ x_k \in \mathbb{X} :  x_k ~\mathrm{is}~\mathrm{feasible}  ~\mathrm{for}~ (7) \big\}. 
\end{equation}
Moreover, we give another two definitions for the $\mathbb{F}$ in \eqref{definition:F}. 
\begin{defn}[Full Region]  \label{Full:region}
	The feasible region $\mathbb{F}$ is full if $\mathbb{F} = \mathbb{X} \subset \mathbb{R}^n$. 
\end{defn}
\begin{defn}[Globally Full Region]  \label{Globally: Full:region}
	The feasible region $\mathbb{F}$ is globally full if $\mathbb{F} = \mathbb{X} = \mathbb{R}^n$. 
\end{defn}
Ideally we would want that the feasible region $\mathbb{F}$ is a globally full region. But it is impossible for both theoretical and application aspects. The best case is that the feasible region $\mathbb{F}$ is a full region such that, we can execute the receding horizon NMPC all over the state space $\mathbb{X}$ with guarantee of the closed-loop stability. Therefore, the main objective of this paper is to learn the feasible region of a nonlinear system under the MPC framework, as well as the nonlinear perturbed case. 
Our method provides freedom in the design of the feasible region needed for closed-loop stability and does not need accurate information of the additively perturbed nonlinear system. 
Furthermore, we give some considerations for enlarging the feasible region so as to make it to be a full region.

\section{Basic Results}
In this section, we use the support vector machine learning method to estimate the feasible region for the nonlinear system given by \eqref{eq:system:nominal} based on the feasibility information of a low-discrepancy sequence. Specifically, using deterministic sampling data, the SVM learning method is employed to learn the boundary function of the feasible region.  

To that end, we let $\partial \mathbb{F}$ denote the boundary of the feasible region $\mathbb{F}$ given by \eqref{definition:F}. Then we make the following assumption. 

\begin{assum}\label{zero:superlevel:set} 
	There exists a continuous function $\phi(\cdot) \in \mathcal{C}^1 (\mathbb{X})$ such that, the feasible region $\mathbb{F}$ can be represented as the zero superlevel-set of $\phi(\cdot)$; that is
	\begin{align} \label{zero:superlevel:set1}
		\mathrm{Int}(\mathbb{F}) = \bigl\{ x \in \mathbb{X}: \phi(x) >0 \bigr\}. 
	\end{align}
\end{assum}
Now using Assumption~\ref{zero:superlevel:set}, to characterize the feasible region $\mathbb{F}$, we just need to construct the function $\phi(\cdot)$. 

\subsection{Point-wise sampling data and SVM learning}
To construct the function $\phi(\cdot)$, we need to sample the state space $\mathbb{X}$ and obtain the point-wise sampling data. To that end, we borrow the standard notion of low-discrepancy sequences from \cite{Niederreiter1988low} and construct a data set on a multilevel sparse-grid. 
\begin{defn} [\cite{chakrabarty2016support}]
	The discrepancy of a sequence $\{x_{k}\}_{k=1}^{N} \subset \mathbb{X}$ is defined as
	\begin{align*}
	\mathscr{D}_{N}\left(  \{x_{k}\}_{k=1}^{N} \right)  \triangleq   \underset{X\in \mathscr{T}} {\mathrm{sup}} \bigg|\frac{ \mathrm{Card}_{\mathbb{X}} (\mathbb{S})}{N} - \frac{\mathrm{Vol}(\mathbb{S})}{\mathrm{Vol}(\mathbb{X})} \bigg|, 
	\end{align*}
	where $\mathscr{T}$ is a set of state dimensional intervals and defined as $\mathscr{T}  \buildrel \vartriangle \over  = \prod_{k=1}^{M} \bigl[a_k,b_k\bigr)=\bigl\{ x_k \in \mathbb{X}:a_k\le x_k<b_k \bigr\} \subset\mathbb{X}$, and $M = \mathrm{Card}_{\mathbb{X}} (\mathbb{S}) \triangleq \mathrm{Card} \bigl\{ k\in\left\{ 1,\ldots,N\right\} : x_k \in \mathbb{S} \bigr\}$. For a low-discrepancy sequence $\{x_{k}\}_{k=1}^{N}$, $\underset{N\to\infty}{\mathrm{lim}} \mathscr{D}_{N} \bigl( \left\{ x_k\right\} _{k=1}^{N} \bigr)  =  0$. 
\end{defn}

Over the state space $\mathbb{X}$, we generate  a point-wise data set
\begin{align}\label{data:set}
\varPi   \buildrel \vartriangle \over  = \left\{ x_k\right\} _{k=1}^{N}  = \bigl\{ x_1, x_2, \ldots, x_N \bigr\},
\end{align}
and assume that it satisfies the following assumption. 

\begin{assum} \label{low:discrepancy:sequence}
The sequence $\left\{ x_k\right\} _{k=1}^{N}$ is a low-discrepancy sequence on $\mathbb{X}$. 
\end{assum}

Now, using Assumption~\ref{low:discrepancy:sequence}, we describe the training data in $\varPi$ following the results presented in \cite{chakrabarty2016support}. 
%
For each point $x_k$, $k=1, \ldots, N$ in $\varPi$, we solve the optimization problem \eqref{eq:7} over a finite-time horizon $T$. 
If there is a feasible solution for $x_k$, then we label it as $\mathcal{O}(x_k) = +1$; otherwise, we label it as $\mathcal{O}(x_k) = -1$. Then we have the following two cases:
\begin{align}\label{Bi:classifier}
\mathcal{O}(x_k)  = \begin{cases}
+1, \quad x_k \in \mathbb{F}; \\
-1, \quad x_k \notin \mathbb{F}, \\
\end{cases}
\end{align}
for each point $x_k$ in $\varPi$.  Note that this is a two-class pattern recognition problem, which divides the data set $\varPi$ into two classes, that is, 
\begin{align*}
\varPi^{+} & =\big\{ x_{k} \in \varPi: ~ \mathcal{O}(x_k)  = +1 \big\};  \\
\varPi^{-} & =\big\{ x_{k} \in \varPi: ~ \mathcal{O}(x_k)  = -1 \big\}. 
\end{align*}
This problem can be solved using the SVM bi-classifier. Here, we use the kernel-based SVM classification method, see e.g. \cite{vapnik2013nature,burges1998tutorial}. Now, we can construct a decision function $\phi_N(\cdot) \in \mathcal{C}^1 (\mathbb{X})$, 
\begin{align}\label{decison:function}
\phi_N(x): \varPi \to y_k \in \{+1, -1\}, 
\end{align}
 that accurately classifies an arbitrary state in $\varPi$ as feasible or infeasible. 

%
%

Next, to facilitate the classification, we introduce a map: $\Gamma: \mathbb{X} \to \mathbb{H}$, which maps the data in $\varPi$ to a higher-dimensional Hilbert space $\mathbb{H}$ where the data is linearly separable, see e.g. \cite{burges1998tutorial}. In this case, we can let $\phi_N(x) = \rho \Gamma(x) + v$, where $\rho$ and $v$ are vector parameters that determine the orientation of the separating hyperplane. Referring to \cite{ong2006enlarging}, our problem can be formulated as the constrained optimization problem 
\begin{subequations}\label{non:constrained:OP}
	\begin{align}\label{non:constrained:OP:target}
	\rho^\star = \underset{ \rho } { \arg\min }  ~  \frac{1}{2} \rho^\mathrm{T} \rho + L \sum_{k =1}^{N} s_k, 
	\end{align}
	subject to
	\begin{align}
	&  \rho  \Gamma(x_k) + v \ge 1 - s_k, \quad  s_k \ge 0, \quad \forall k \in \varPi^{+},   \label{non:constrained:OP:1}  \\
	&  \rho  \Gamma(x_k) + v \le -1 + s_k, \quad \forall k \in \varPi^{-},   \label{non:constrained:OP:2} 
	\end{align}
\end{subequations}
where $L$ is a regularization parameter and $s_k$ is a slack variable to relax separability constraints. 

The optimization problem \eqref{non:constrained:OP} is convex and can be solved via its dual formulation. Letting $\alpha_i$ and $\alpha_j$ denote the Lagrange multiplier vector for \eqref{non:constrained:OP:1} and \eqref{non:constrained:OP:2}, \eqref{non:constrained:OP} can reformulated as the dual problem
\begin{subequations}\label{non:constrained:OP2}
	\begin{align}\label{non:constrained:OP2:target}
\alpha^{\star} = \underset{\alpha}{\mathrm{arg~min}}  &  ~\frac{1}{2}\sum_{i=1}^{N}\sum_{k=1}^{N}\alpha_i \alpha_{k} y_{i}y_{k} \Gamma(x_i) \Gamma(x_k) - \sum_{k=1}^{N}\alpha_k, \\
\mathrm{subject\;to}~~ 
	&  \sum_{k=1}^{N}\alpha_k y_{k}=0, \label{non:constrained:OP2:1}  \\
	&   0 \le \alpha_{k} \le L,  \quad  k = 1, \ldots, N. \label{non:constrained:OP2:2}
	\end{align}
\end{subequations}
After solving \eqref{non:constrained:OP2}, we obtain the SVM decision function as
\begin{equation}\label{SVM:function}
\phi_N(x)  = \sum_{k=1}^{N}\alpha_{k}^{\star} y_{k} \Gamma(x_k) \Gamma(x) + v, 
\end{equation}
and the estimated feasibility region boundary is given by
\begin{equation} \label{SVM:functionregion}
\mathcal{S}_N  \triangleq   \left\{ x\in\mathbb{X}:  \phi_N(x) = 0 \right\}. 
\end{equation}

To improve accuracy of the SVM learning when effective data in $\varPi^{+}$ are not enough, we iteratively select new points satisfying $\mathcal{O}(x_k)  = +1$ and expand into the data set $\varPi^{+}$ to ensure richness of data. 
This ensures that a non-trivial decision function \eqref{SVM:function} exists. 
Furthermore, for the inner product $\Gamma(x_k) \Gamma(x)$ in \eqref{SVM:function}, it can be replaced by using the kernel function. A common choice of kernel function is the Gaussian kernel: 
	\[
\Gamma(x_k) \Gamma(x) = \frac{{\|x_{k}-x\|}^2}{2\sigma^{2}}, 
\]
where $\sigma^2$ is the kernel variance, see e.g. \cite{steinwart2001influence}.

\subsection{Inner and outer approximations of $\partial\mathbb{F}$}
The outer and inner approximations of a bounded region were introduced in \cite{darup2012low,deffuant2007approximating}. Using \eqref{zero:superlevel:set1}, we define the 
strict outer ($\mathbb{F}^{-}$) and inner ($\mathbb{F}^{+}$) region of the feasible region $\mathbb{F}$ as
\begin{align}
\mathbb{F}^{+} & \triangleq \left\{ x\in\mathbb{X}:\phi(x)>\varepsilon\right\};  \label{F+}\\ 
\mathbb{F}^{-}  & \triangleq  \left\{ x\in\mathbb{X}:\phi(x)< - \varepsilon\right\}.  \label{F-} 
\end{align}
for some $\varepsilon >0$. 
In this subsection, we construct strict outer and inner approximations of the feasible region $\mathbb{F}$, which ensures that they contain no feasible (infeasible) samples. 

Now, based on the super-level sets of $\phi_N(x)$, we can find a positive constant $\varepsilon^{+}>0$ which defines the inner approximation $\mathbb{F}_N^{+}$ by 
\begin{align}
\mathbb{F}_N^{+} & \triangleq \bigl\{ x\in \mathbb{X}: \phi_N(x) > \varepsilon^{+}\bigr\}. \label{eq:outer}
\end{align}
Using the data in $\varPi^+$, we determine the parameter $\varepsilon^{+}$ by solving the optimization problem
\begin{subequations}\label{eq:outer:OP}
	\begin{align}
	\varepsilon^{+}  & =\underset{\varepsilon>0}   {\mathrm{arg~min}} ~ \varepsilon,  \label{eq:outer:OP:1}  \\
	\mathrm{subject\;to}~ & \phi_N (x_{k} ) > \varepsilon,  \quad\forall x_{k}\in \varPi^+. \label{eq:outer:OP:2}
	\end{align}
\end{subequations}

Analogously, there exists a negative constant $\varepsilon^{-} <0$ which defines the outer approximation $\mathbb{F}_N^{-}$ by 
\begin{align}
\mathbb{F}_N^{-} & \triangleq \bigl\{ x\in  \mathbb{X}: \phi_N(x) < \varepsilon^{-}\bigr\}.  \label{eq:inner}
\end{align}
This constant $\varepsilon^{-}$ can be determined by using the data in $\varPi^-$ and solving the optimization problem 
\begin{subequations}\label{eq:inner:OP}
	\begin{align}
	\varepsilon^{-}  & =\underset{\varepsilon<0}   {\mathrm{arg~max}} ~ \varepsilon,  \label{eq:inner:OP:1}  \\
	\mathrm{subject\;to}~ & \phi_N (x_{k} ) < \varepsilon,   \quad\forall x_{k}  \in \varPi^-.  \label{eq:inner:OP:2}
	\end{align}
\end{subequations}

We highlight that novel algorithms for constructing strict inner and outer approximations of the feasible region $\mathbb{F}$ are developed by solving \eqref{eq:outer:OP} and \eqref{eq:inner:OP}. Besides, the inner approximation $\mathbb{F}_N^{+}$ is based only on the feasible samples in $\varPi^+$ and the outer approximation $\mathbb{F}_N^{-}$ is based only on the infeasible samples in $\varPi^-$. 

\subsection{Convergence analysis}
%
In this subsection, we provide convergence analysis for the approximation error. 
We show that as the number of samples $N$ from the low-discrepancy sequence increases, the inner and outer approximations with universal kernels will converge to a strict approximation of the actual feasible region $\mathbb{F}$. 

First, we introduce definition of universal kernels adopted from \cite{steinwart2001influence} and give a proposition. 
\begin{defn} \label{universal:kernel}
	A continuous kernel $\mathcal{K}(\cdot, \cdot)$ is universal 
	if the space of all functions induced by $\mathcal{K}(\cdot, \cdot)$ is dense in $\mathcal{C}(\mathbb{X})$. 
\end{defn}
\begin{prop} \label{universal:kernel:property}
	Every universal kernel separates all compact subsets of $\mathbb{X}$.
\end{prop}
Definition~\ref{universal:kernel} indicates that for each continuous function $\phi$ and $\varepsilon>0$, 
there exists a function $\phi_{\mathcal{K}}$ induced by a universal kernel $\mathcal{K}$ such that $|\phi-\phi_{\mathcal{K}}|<\varepsilon$.
Proposition~\ref{universal:kernel:property} indicates that the kernel-based SVM classification method can be adopted to solve the feasibility boundary learning problem. 
In this case, if we design 
\begin{align} \label{kernel:function}
	\mathcal{K}(x_k, x) = \Gamma(x_k) \Gamma(x),
\end{align}
then we can accurately separate the pairwise disjoint compact sets $\varPi^+$ and $\varPi^-$ in $\mathbb{X}$. 

We highlight that there is not computational method for determining such a separating function in \cite{steinwart2001influence}. 
By modifying the arguments presented in \cite{steinwart2001influence}, we collect the training data from a low-discrepancy sequence. 
It guarantees that increasing the number of samples $N$ for
the SVM with a universal kernel results in an increasing accuracy of the feasible region learning, and reduces the conservativeness of the inner and outer approximations of $\mathbb{F}$. 
We summarize our discussions as the following theorem, which is slightly similar to the results given in \cite{chakrabarty2016support}. 

\begin{thm}\label{theprem:main}
	Consider the system \eqref{eq:system:nominal}, suppose Assumptions~\ref{eq:J:l}, \ref{eq:J:II}, \ref{eq:J:F}, \ref{zero:superlevel:set}, and \ref{low:discrepancy:sequence} hold, and let $\mathcal{K}(\cdot, \cdot)$ be a universal kernel on $\mathbb{X}$ in \eqref{kernel:function} and $\phi_{N, \mathcal{K}} (x)$ be the kernel-based SVM decision function from \eqref{SVM:function}. Then for any $L>0$ in the optimization problem \eqref{non:constrained:OP2} and for two compact subsets $\mathbb{F}^{+} \subset \mathbb{X}$ and $\mathbb{F}^{-} \subset \mathbb{X}$ in \eqref{F+} and \eqref{F-}, 
	\begin{align}
	& \underset{N \to \infty}{ \mathrm{lim}} \phi_N(x) = \phi(x), \\
	&	\underset{N \to \infty}{ \mathrm{lim}}  \mathbb{F}_N^+ = \mathbb{F}^+, \quad \underset{N \to \infty}{ \mathrm{lim}}  \mathbb{F}_N^- = 	\mathbb{F}^-, \\
	&	\underset{N \to \infty}{ \mathrm{lim}} \mathcal{S}_N  = \partial \mathbb{F}. 
	\end{align}
\end{thm}

Theorem~\ref{theprem:main} indicates that there exists a sufficiently large number of samples $N_{0}$ such that if $N\ge N_{0}$, 
\[
\begin{cases}
\phi(x) > 0,  &  \quad  \mathrm{for} \; x\in\mathbb{F}^{+};  \\
\phi(x) < 0,  &  \quad  \mathrm{for} \; x\in\mathbb{F}^{-}, 
\end{cases}
\]
and the strict inner and outer approximations of the feasible region $\mathbb{F}$ can be obtained by sampling $\mathbb{X}$ using low-discrepancy sequences. 
Therefore, the SVM classifier using universal kernels is able to approximate our feasible region boundaries with sufficiently high accuracy. 

\section{Main Results}
In Section~3, we present a method to estimate the feasible region where the states of the nominal system \eqref{eq:system:chen1997} can be driven to the terminal set $\Omega$ while satisfying the constraints \eqref{eq:system:description:2}. In this section, we extend this result to the perturbed nonlinear system given by \eqref{eq:system:description}. However, the region $\mathbb{F}$ obtaining from Section 3 may not be feasible for the perturbed system when there exists unknown disturbance. 

Using Assumption A2, \eqref{eq:system:description}, and \eqref{eq:system:nominal}, it follows that 
\begin{align}
| F(x(t), u(t), w(t)) - f(x(t), u(t)) | \le \bar{w}. 
\end{align}
So for the perturbed system \eqref{eq:system:description}, the set $\Omega$ may not be positively invariant with the local feedback law $u(t)$. In this case, we introduce the robust control invariant (RCI) set for \eqref{eq:system:description} following the results presented in \cite{sun2018robust}. 
\begin{defn}[RCI set]  \label{RCIset}
	Consider the perturbed system \eqref{eq:system:description} with constraints $x \in \mathbb{X}$, $u\in \mathbb{U}$ and $w \in \mathbb{W}$. The set $\bar{\Omega} \subset \mathbb{X} \subseteq \mathbb{R}^n$  is an RCI set, if there exists an ancillary controller $\kappa(x)$ with $\bar{u} = u + \kappa(x) \in \mathbb{U}$ such that $\forall \bar{x} \in \bar{\Omega}$ and $\forall w \in \mathbb{W}$, $x(t) \in \bar{\Omega}$. 
\end{defn}

For the design of the RCI set $\bar{\Omega}$ (terminal region) for \eqref{eq:system:description}, the set-theoretic methods  and tube-based MPC scheme are very useful, see e.g. \cite{blanchini:2007} and \cite{Rakovic:2012}. However, the method developed in \cite{Rakovic:2012} can not be applied to our problem, since the perturbed system given by \eqref{eq:system:description} is continuous-time and nonlinear case; the authors in \cite{ong2006enlarging} and  \cite{chakrabarty2016support}  only consider the nonlinear systems, not for the perturbed system. 

\subsection{RCI set design based on set-theoretic method}

Recalling that the disturbance in \eqref{eq:5} is actually an additive disturbance and satisfies
\begin{align} \label{eq:w}
w (t) = F(x(t), u(t), w(t)) - f(x(t), u(t)), ~ t \in \mathbb{ R}_{\ge 0}. 
\end{align}
Besides, we note that the state and control do not depend on the disturbance $w$. 

Now, using Assumption~A2, we introduce the following set
\begin{align} \label{eq:Set}
	\Omega_0 \triangleq \bigl\{ x \in \Omega: |x|_{\partial  \Omega  } \ge \bar{w} \bigr\}, 
\end{align}
where $\Omega$ is a terminal region for \eqref{eq:system:nominal} and $\partial  \Omega $ is the boundary of the set $\Omega$.  
Then for \eqref{eq:Set}, we have the following lemma.

\begin{lem}  \label{Lem:Omega}
	Consider the set \eqref{eq:Set} and assume that $\Omega \supset \mathbb{R}^n (\bar{w})$ holds. Then, $\Omega_0$ is a RCI set for the perturbed system \eqref{eq:system:description}. 
\end{lem}
\begin{pf}
	Note that $\Omega$ is a control invariant set for the nonlinear system \eqref{eq:system:nominal} and $\Omega_0 \subset \Omega$. Since $\Omega \supset \mathbb{R}^n (\bar{w})$, $\Omega \setminus \Omega_0 \neq \varnothing$. Then using \eqref{eq:Set}, we have that for any $x \in \Omega_0$, $ |x + w |_{ \partial  \Omega_0 } \ge |x|_{ \partial  \Omega } - \bar{w} \ge 0 $ for all $ w \in \mathbb{W}$. That is, for any $x \in \Omega_0 $, $x + w \in \Omega$ for all $ w \in \mathbb{W}$. According to Definition~\ref{RCIset}, $\Omega_0$ is a RCI set for \eqref{eq:system:description}. 
\end{pf}

Lemma~\ref{Lem:Omega} indicates that the RCI set $\bar{\Omega}$ for \eqref{eq:system:description} can be obtained using the terminal set $\Omega$ for \eqref{eq:system:nominal}, while the set $\Omega$ can be computed using the results developed in \cite{chen1998quasi} for perturbed nonlinear  systems.

\subsection{Feasible region estimation for robust NMPC}
In this subsection, we develop a feasible region estimation method for \eqref{eq:system:description} under the NMPC scheme using the feasibility information of the low-discrepancy data samples. 

First, using the RCI set $\Omega_0$ in \eqref{eq:Set} for the perturbed system given by \eqref{eq:system:description}, we modify the problem \eqref{eq:7} to a new optimization problem of robust NMPC which is given by
	\begin{subequations}\label{eq:robust:NMPC}
	\begin{align}\label{robust:NMPC}
	\hat{u}^\star = \underset{ u } { \arg\min }  ~   \mathcal{J} (x_k,  \hat{u}(\cdot)), 
	\end{align}
	subject to
	\begin{align}
	&  \dot{\hat x} (\tau) =  f(\hat{x} (\tau), \hat{u} (\tau)),  \quad 0 \le \tau \le T,   \label{eq:robust:NMPC:A}  \\
	&   \hat{x} (\tau) \in \mathbb{X}, \quad \hat{u} (\tau) \in \mathbb{U}, \quad 0 \le \tau \le T,   \label{eq:robust:NMPC:B} \\
	&  \hat{x} (0) = x_k,   \quad \hat x (T) \in \Omega_0,    \label{eq:robust:NMPC:C}
	\end{align}
\end{subequations}
where $\mathcal{J}(x_k,  \hat{u}(\cdot))$ is given by \eqref{Perf} and $\Psi(\cdot, \cdot)$ and $\Phi(\cdot)$ satisfy Assumptions~\ref{eq:J:l} and \ref{eq:J:F}. 

Next, using the same procedure developed in Section 3, we generate the data set and construct an approximate function of the boundary for the feasible region. Then similar to Theorem~\ref{theprem:main}, we have the following theorem which shows the approximation of the feasible region $\bar{\mathbb{F}}$ for the perturbed systems in \eqref{eq:system:description}. 

\begin{thm}\label{theprem:main:two}
		Consider the system \eqref{eq:system:description}, suppose Assumptions~\ref{eq:J:l}, \ref{eq:J:F}, \ref{zero:superlevel:set}, and \ref{low:discrepancy:sequence} hold, and let $\bar{\Omega}$ be the maximal RCI set in $\mathbb{X}$
	and $\bar{\mathbb{F}}$ be the maximal feasible region for \eqref{eq:system:description}. Then, 
	$\Omega_0 \subset \bar{\Omega}$ and 
	\begin{align} \label{eq:Region}
    \bar{\mathbb{F}}_0 \triangleq \bigl\{ x \in {\mathbb{F}}: |x|_{\partial  {\mathbb{F}}  } \ge \bar{w} \bigr\} \subset \bar{\mathbb{F}}. 
	\end{align}
\end{thm}
Theorem~\ref{theprem:main:two} indicates that based on the feasible region $\mathbb{F}$ for \eqref{eq:system:nominal}, we can obtain the conservative terminal region $\Omega_0$ and the conservative feasible region $\bar{\mathbb{F}}_0$ for the perturbed systems. Since the disturbance $w$ in \eqref{eq:w} is unknown, it is very difficult to compute the strict feasible region for the perturbed system \eqref{eq:system:nominal}. Therefore,  we have adopted the low-discrepancy sequences, optimal kernel-based SVM classification, and  set-theoretic methods to increase the accuracy of the data-driven feasible region estimation. Although the feasible region $\bar{\mathbb{F}}$ can be enlarged with a longer prediction horizon, it needs more computation resources for solving the optimization problems \eqref{eq:robust:NMPC}.


\section{illustrative numerical example}\label{sec:examples}
To illustrate the key ideas presented in this paper, we consider a perturbed nonlinear system representing the cart-and-spring system adopted from \cite{magni2003robust,li2014distributed} given by 
\begin{equation*}\begin{cases}
\dot{x}_{1} (t) =  x_{2}(t), \\
\dot{x}_{2} (t) =  - \frac{k_0}{M} e^{- x_1(t)} x_1(t)  - \frac{h_d}{M} x_2(t) + \frac{u(t)}{M} + w(t), 
\end{cases}
\end{equation*}
where $x = \left[x_1; x_2 \right] \in \mathbb{R}^2$, $u \in \mathbb{R}$,  and $w \in \mathbb{R}$. The parameters are given as: $M = 1.8~\mathrm{kg}$, $k_0 = 1.2~\mathrm{N/m}$, $h_d =0.25~\mathrm{N\cdot s/m}$,  
%
The constrained state and control space associated with the system is given by $\mathbb{X} = \{ x \in \mathbb{R}^2: |x| \le 2 \}$ and $\mathbb{U} = \{ u \in \mathbb{R}: |u| \le 3\}$ and let $w(t) \in \mathbb{B}^1(0.01)$. 
	To satisfy Assumptions~\ref{eq:J:l} and \ref{eq:J:F}, we set the stage cost as $\Psi(x, u) = \|x\|_{Q}^2 + \|u\|_R^2$ with $Q =I_2$ and $R = 0.1$ and the terminal cost as $\Phi(x) = \|x\|_{P}^2$ with $P = \left[5.0511, -2.2731; -2.2731, 2.4586\right]$. The local controller is designed as $\kappa(x) = -Kx$ with $K = \left[4.2291~4.8551\right]$ using the procedure presented in \cite{chen1998quasi}.  
	
\begin{figure}[ht]
	\includegraphics[width=.99\linewidth]{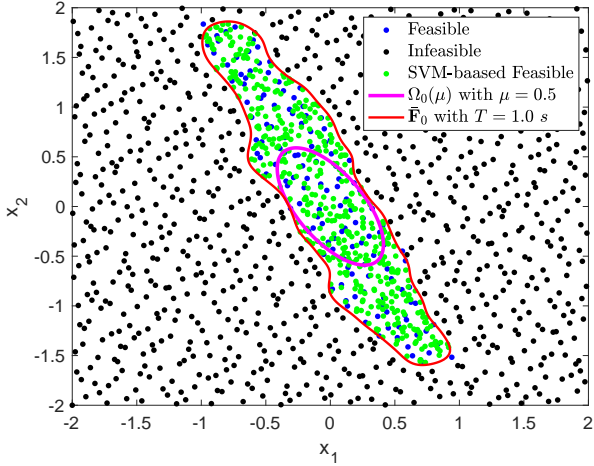}
	\caption{The boundary of the feasible region~$\bar{\mathbb{F}}_0$ resulting from the optimization  \mbox{(\ref{eq:robust:NMPC})} with the terminal region $\Omega_0(0.5)$ and the prediction horizon $T = 1.0~s$ in \textcolor{red}{red}. 
		The infeasible samples in $\varPi^{-}$ are shown in \textcolor{black}{black}; 
		the feasible samples in $\varPi^{+}$ are shown in \textcolor{blue}{blue}; 
		the feasible samples selected using \eqref{SVM:functionregion} are shown in \textcolor{green}{green}; 
	   the boundary of the terminal region $\Omega_0(0.5)$ is shown in \textcolor{magenta}{magenta}. 
   }
	\label{figure:feasible:region:SVM:T10s}
	\vspace{4mm}\end{figure}

\begin{figure}[ht]
	\includegraphics[width=.99\linewidth]{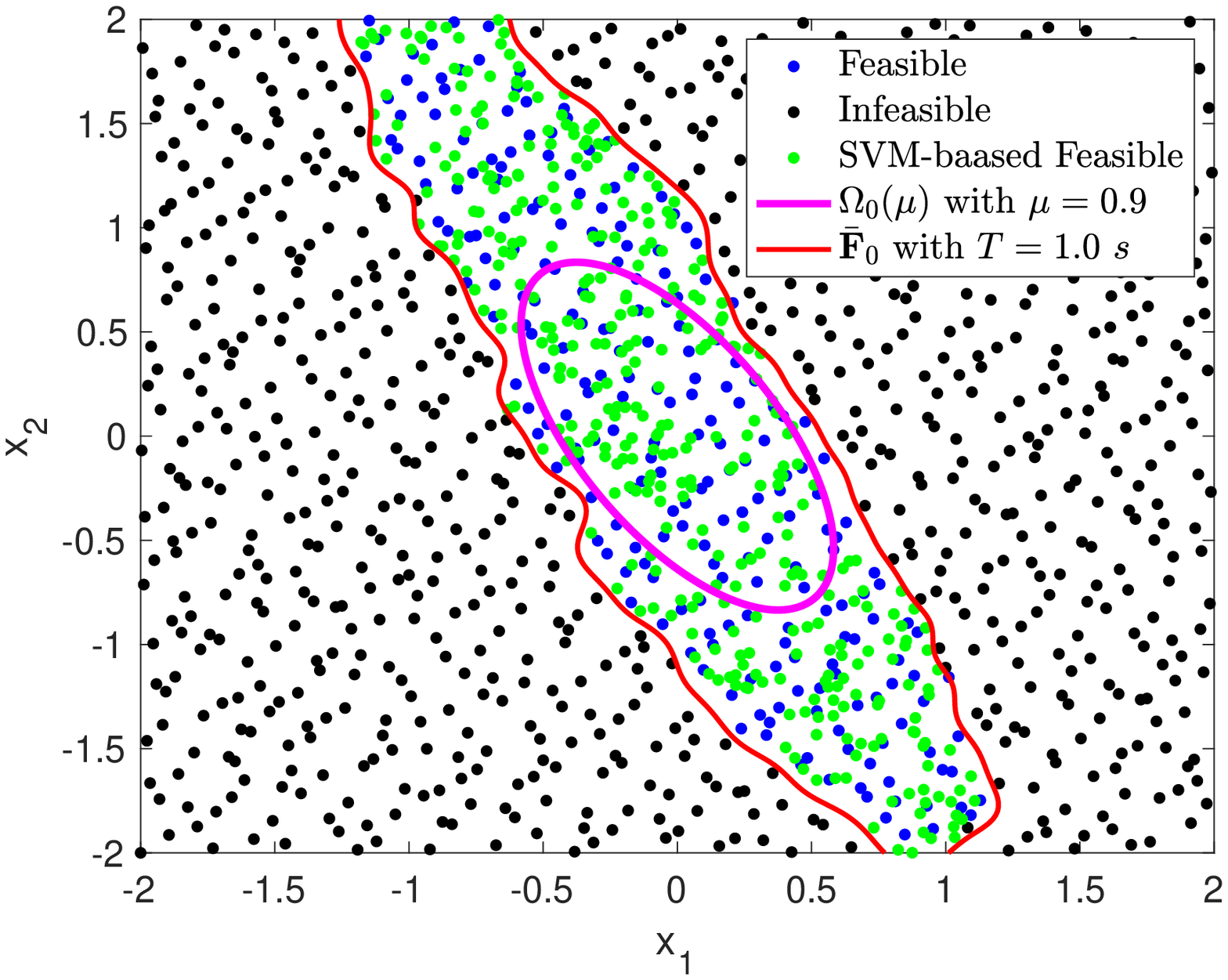}
	\caption{The boundary of the feasible region~$\bar{\mathbb{F}}_0$ resulting from the optimization  \mbox{(\ref{eq:robust:NMPC})} with the terminal region $\Omega_0(0.9)$ and the prediction horizon $T = 1.0~s$ in \textcolor{red}{red}. 
	The other notations are the same as Fig.~\ref{figure:feasible:region:SVM:T10s}. 
}
	\label{figure:feasible:region:SVM:T15s}
\end{figure}

\begin{figure}
\includegraphics[width=0.99\linewidth]{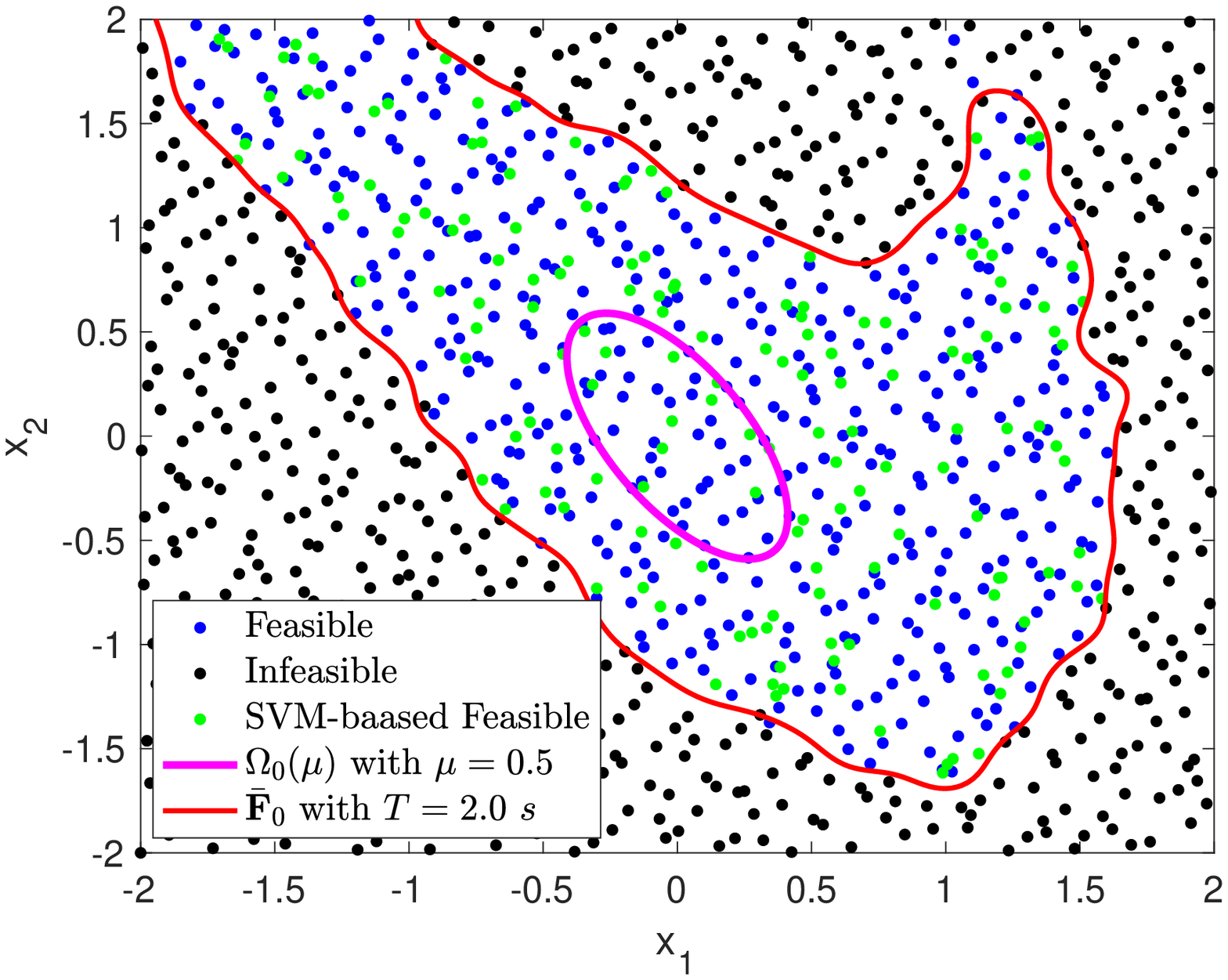}
\caption{The boundary of the feasible region~$\bar{\mathbb{F}}_0$ resulting from the optimization  \mbox{(\ref{eq:robust:NMPC})} with the terminal region $\Omega_0(0.5)$ and the prediction horizon $T = 2.0~s$ in \textcolor{red}{red}. 
The other notations are the same as Fig.~\ref{figure:feasible:region:SVM:T10s}. 
}
\label{figure:feasible:region:SVM:T20s}
\end{figure}

	Next, the terminal region $\Omega$ is parametrized by $\Omega(\mu) = \bigl\{ x \in \mathbb{R}^2: x^\mathrm{T} P x \le \mu \bigr\}$, where $\mu \in (0, 1)$. 
	Then, we use a kernel-based SVM with $\sigma = 0.8$ for the feasible region approximation. Using the feasibility information of data samples, approximations of the feasible  region boundaries with different prediction horizons $T$ and different $\mu$  are illustrated in Figures.~\ref{figure:feasible:region:SVM:T10s}--\ref{figure:feasible:region:SVM:T20s}. 
	Furthermore, it can be seen from Figures~\ref{figure:feasible:region:SVM:T10s} and \ref{figure:feasible:region:SVM:T15s} that larger terminal region results in larger feasible region. From Figures~\ref{figure:feasible:region:SVM:T10s} and \ref{figure:feasible:region:SVM:T20s}, we can see that longer prediction horizon results in larger feasible region. 


\section{Conclusions and future work}\label{sec:conclusion}
This paper develops a data-driven learning framework for designing the terminal region and approximating the
feasible region offline using the feasibility information of low-discrepancy data samples and support vector machine
learning. Our approach provides the freedom in the design of the feasible region needed for closed-loop stability, when
disturbance is subject to deterministic bounds. Future work will combine the offline feasible region learning with
online NMPC design. 

%

\bibliography{MyBibliography}
\end{document}